\documentstyle[11pt,amssymb]{article}

\newcommand{\realtilde}{{\char'176}}
\newcommand{\lbrk}{{\linebreak[0]}}
\newcommand{\ds}{\displaystyle}
\newcommand{\op}{{\rm op}}
\newcommand{\erf}{{\rm erf}}
\newcommand{\ts}{\textstyle}
\newcommand{\atan}{{\rm atan}}
\newcommand{\ul}{\underline}

\begin{document}

\begin{center}

{\large\bf UNRESTRICTED ALGORITHMS FOR\\[2ex]
ELEMENTARY AND SPECIAL FUNCTIONS\\[2ex]
Invited Paper}
\footnote{First appeared in
{\em Information Processing 80}
(edited by S.~H.~Lavington),
North-Holland, Amsterdam, 1980, 613--619.
Retyped with minor corrections by
Frances Page at Oxford University Computing
Laboratory, 1999.\\
Copyright \copyright\ 1980--2010, R.~P.~Brent.
\hspace*{\fill} rpb052 typeset using \LaTeX.}\\[6ex]

{\bf Richard P Brent}\\[6ex]
Department of Computer Science\\
 Australian National University\\
Canberra, Australia\\ \vspace*{6mm}
\end{center}

\begin{quote}
We describe some ``unrestricted'' algorithms which are useful for the
computation of elementary and special functions when the precision
required is not known in advance.  Several general classes of
algorithms are identified and illustrated by examples.  Applications
of such algorithms are mentioned.
\end{quote}

\section{Introduction}
\thispagestyle{empty}

Floating-point computations are usually performed with fixed precision:
the machine used may have ``single'' or ``double'' precision 
floating-point hardware, or on small machines fixed-precision
floating-point operations may be implemented by software or firmware.
Most high-level languages support  only  a small number of
floating-point precisions, and those which support an arbitrary number
usually demand that the precision be determinable at compile time.\\

We say that an algorithm has precision $n$ if its result is computed
with error $O(2^{-n})$.  Usually we are interested in the relative
error, but in some cases (e.g.~the computation of $\sin(x)$ for 
$x \simeq \pi$) it is more appropriate to consider the absolute
error.\\

In certain applications it is desirable that the precision of
floating-point operations should be able to be varied at runtime.
In this paper we consider algorithms which may be used to evaluate
elementary and special functions to precision $n$, where $n$ may
be arbitrarily large.  Such algorithms have been termed 
``unrestricted'' by Clenshaw and Olver [13].  Note that algorithms
which are ``unrestricted'' in our sense may have domain restrictions
(e.g.~an ``unrestricted'' algorithm for $\exp(x)$ might be applicable
only for $x \ge 0$), although such restrictions can often be
circumvented by the methods of Section~4, or by combining several
algorithms with different domain restrictions.\\

Unrestricted algorithms depend on the availability of
variable-precision floating-point arithmetic.  At present this is
usually implemented by software, e.g. in the MP package [7],
but it could be implemented in firmware or hardware.  (Note the
historical example of the IBM 1620.)\\

In the following sections we ignore the possibility of floating-point 
underflow or overflow.  For ``ideal'' variable-precision arithmetic
the exponent range should tend to infinity with the\linebreak

\noindent precision~$n$.
For most purposes, though, a whole-word exponent, as used in MP
[7], is adequate to avoid overflow problems.\\

Applications of variable-precision floating-point arithmetic include:

\begin{enumerate}

\item generation and testing of accurate tables of constants
(e.g.~coefficients in minimax polynomial or rational approximations
[9, 25]);

\item computation using numerically unstable algorithms [2, 9];

\item interval arithmetic, where the final intervals may be too
large if fixed-precision arithmetic is used [20, 28, 37];

\item truly machine-independent floating-point computations;

\item testing of floating-point hardware for correctness and
conformity to standards, e.g.~those proposed in [14, 29].

\item number-theoretic computations where very high precision may
be essential [24, 30].

\end{enumerate}

In this paper we concentrate on unrestricted algorithms rather than
their applications.\linebreak
Section~2 summarises some preliminary results.
In each of Sections 3 to 9 we illustrate, by one or two simple
examples, a useful general method leading to unrestricted algorithms.
Some more specialised methods are mentioned in Section~10.  The
field is vast and we make no attempt to be comprehensive.  For
simplicity we usually restrict our attention to real variables
and omit details of the rounding error analysis.  We  also omit any
discussion of desirable high-level language facilities to support
variable precision arithmetic, for which see [10, 19, 32].\\

Very few of the algorithms given below are new, in fact most of the
identities underlying them may be found in [1] or [35].  What
may be new is our viewpoint.  Often an excellent unrestricted 
algorithm is unsuitable for fixed-precision computation, and {\em vice
versa\/}.

\section{Basic arithmetic operations}

We assume that variable-precision floating-point numbers are represented
by an integer exponent and a fraction with $t$ digits to base
$\beta > 1$.  We call such numbers ``precision $n$'' numbers if
$n \simeq (t-1)\log_2\; \beta$.  Addition and subtraction of such
numbers is straight-forward, and requires $O(n)$ operations
[22, 23].  We assume at least one guard digit, so the relative
error in the computed result is at most $\beta^{1-t} = O(2^{-n})$
(see [36]).\\

Let $M(n)$ be the number of operations required for multiplication
of precision $n$ numbers.  By the Sch\"{o}nhage-Strassen algorithm
[22, 33]
$$
M(n) = O(n \log n \log\log n) \;.
\eqno{(1)}$$

For the moderate values of $n$ which usually arise in applications,
an efficient implementation of the classical $O(n^2)$ algorithm may
be faster than the Sch\"{o}nhage-Strassen algorithm or other
asymptotically fast algorithms.\\

Let $D(n)$ be the number of operations required for division of 
precision $n$ numbers.  Under plausible assumptions it may be shown
that $D(n) = O(M(n))$ (see, for example, [5]).  In practice\linebreak
\pagebreak %

\noindent the ``schoolboy'' algorithm, 
which requires $O(n^2)$ operations, may be the
fastest unless $n$ is rather large.\\

It is important to distinguish between multiplication of two precision
$n$ numbers and multiplication or division of a precision $n$
number by a small (single-precision) integer.  The latter require
only $O(n)$ operations if implemented in the obvious way.

\section{Power series}

If $f(x)$ is analytic in a neighbourhood of some point $c$, an obvious
method to consider for the evaluation of $f(x)$ is summation of the
Taylor series
$$ 
f(x)= \sum^{k-1}_{j=0}\; (x-c)^j\; f^{(j)}\; (c)/j! + R_k(x,c) \;.
\eqno{(2)}$$

As a simple but instructive example we consider the evaluation of
$\exp(x)$ for $|x|\leq 1$, using
$$ 
\exp(x) = \sum^{k-1}_{j=0}\; x^j/j! + R_k(x) \;,
\eqno{(3)}$$

where $|R_k(x)|\leq e/k!$\\

Using Stirling's approximation for $k!$, we see that $k \ge K(n) \sim
n/\log_2 n$ is sufficient to ensure that $|R_k(x)|=
O(2^{-n})$.  Thus the time required is $O(nM(n)/\ln n)$.\\

In practice it is convenient to sum the series in the forward direction
$(j=0, 1, \ldots, k-1)$.  The terms $T_j = x/j!$ and partial sums
$$S_j = \sum^j_{i=0} T_i$$
may be generated by the recurrence
$T_j = x \times T_{j-1}/j$, $S_j=S_{j-1}+T_j$, and the summation
terminated when $|T_k|< 2^{-n}$.  Thus, it is not necessary
to estimate $k$ in advance, as it would be if the series were
summed by Horner's rule in the backward direction
$(j=k-1, k-2, \ldots, 0)$.\\

We now consider the effect of rounding errors, under the assumption
that floating-point operations satisfy
$$ 
fl(x \;\op\; y) = (x \;\op\; y)(1+\delta) \;,
\eqno{(4)}$$

where $|\delta|\leq \varepsilon$ and ``op'' = ``+'',
``$-$'', ``$\times$'' or ``/''. Here  $\varepsilon \leq \beta^{1-t}$
is the ``machine-precision'' [36].
Let $\widehat{T}_j$ be the computed value of $T_j$, etc.  Thus
$$
|\widehat{T}_j - T_j|\; / \;|T_j|\;\leq\; 2j \varepsilon
+ O(\varepsilon^2)
\eqno{(5)}$$

and
$$
|\widehat{S}_k - S_k|\;\leq\; ke\varepsilon  + \sum^k_{j=1}\;
                                     2j \varepsilon|T_j|+
                                       O(\varepsilon^2)
$$
\vspace*{-9mm}
\begin{flushright}
(6)
\end{flushright}
\vspace*{-4mm}
$$ \hspace*{20mm} \;\leq \;(k+2)e \varepsilon +
                           O(\varepsilon^2) = O(n\varepsilon)\;.
$$

Thus, to get $|\widehat{S}_k - S_k|= O(2^{-n})$ it is sufficient
that $\varepsilon = O(2^{-n}/n)$, i.e.~we need to work with about
$\log_\beta n$ guard digits.  This is not a significant overhead if 
(as we assume) the number of digits may vary dynamically.  The slightly
better error bound obtainable for backward summation is thus of no 
importance.\\

In practice it is inefficient to keep $\varepsilon$ fixed.  We can
profitably reduce the working precision when computing $T_k$ from
$ T_{k-1}$ if $|T_{k-1}|\ll 1$, without significantly
increasing the error bound.\\

It is instructive to consider the effect of relaxing our restriction
that $|x|\leq 1$.  First suppose that $x$ is large and
positive.  Since $|T_j|>|T_{j-1}|$ when
$j <|x|$, it is clear that the number of terms required
in the sum (3) is at least of order $|x|$. Thus, the method is
slow for large $|x|$ (see Section~4 for faster methods in
this case).\\

If $|x|$ is large and $x$ is negative, the situation is even
worse.  From Stirling's approximation we have
$$
\max_{j \ge 0}\;|T_j|\; \simeq\; \frac{\exp|x|}{\sqrt{2\pi|x|}} \;,
\eqno{(7)}$$

but the result is $\exp(-|x|)$, so about
$2|x|/ \ln \beta$ guard digits are required to compensate
for Lehmer's ``catastrophic cancellation'' [15].  Since
$\exp(x) = 1/\exp(-x)$, this problem may easily be avoided,  but
the corresponding problem is not always so easily avoided for other
analytic functions.\\

In the following sections we generally ignore the effect of rounding
errors, but the results obtained above are typical.  For an example of
an extremely detailed error analysis of an unrestricted algorithm,
see [13].\\

To conclude this section we give a less trivial example where power
series expansions are useful.  To compute the error function
$$ 
\erf(x) = 2\pi^{-{1/2}} \int^x_0 e^{-u^{2}}\; du \;, $$

we may use the series
$$ 
\erf(x) = 2\pi^{-{1/2}}\; \sum^\infty_{j=0}\;
    \frac{(-1)^j\; x^{2j+1}}{j!(2j+1)}
\eqno{(8)}$$

or
$$ \erf(x) = 2\pi^{-{1/2}} \exp(-x^2)
  \; \sum^\infty_{j=0} \;\frac{2^j \; x^{2j+1}}
                          {1 \cdot 3 \cdot 5\cdots(2j+1)} \;.
\eqno{(9)}$$

The series (9) is preferable to (8) for moderate $|x|$ because
it involves no cancellation.  For large $|x|$ neither series
is satisfactory, because $\Omega(x^2)$ terms are required, and it is
preferable to use the asymptotic expansion or continued fraction
for $\erf c(x) = 1 - \erf(x)$: see Sections 5 and 6.

\section{Halving identities}

In Section~3 we saw that the power series is not suitable for
evaluation of $\exp(x)$ if $|x|$ is large.  To reduce the size of
the argument we may use the identity
$$ 
\exp(x) = [\exp(x/2)]^2		%
\eqno{(10)}$$

as often as necessary.  When applied $k$ times, (10) gives
$$
\exp(x) = [\exp(2^{-k}x)]^{2^k} \;.
\eqno{(11)}$$

If $k = \lfloor cn^{{1/2}} \rfloor + \log_2 |x|$
for some positive constant $c$, and (11) is used in conjunction with
the power series algorithm of Section~3, the time required to evaluate
$\exp(x)$ to precision $n$ for large $|x|$ is
$$ O[(n^{{1/2}} + \ln|x|)M(n)]\;,$$ better than the
$O[(n/\ln n+|x|)M(n)]$ result of Section~3 (the case
$k=0$).\\

Similar ``halving'' (or ``doubling'') identities, derived by replacing
$x$ by $ix$ in (10),   may be used to evaluate trigonometric and
inverse trigonometric functions [5, 7, 13, 32].\\

Other identities are useful in special applications:  see Section~10
for some examples.

\section{Asymptotic expansions}

Rarely does a single method suffice to evaluate a special function
over its whole domain.  For example, the exponential integral
$$ 
E_1(x) = \int^\infty_x \frac{\exp(-u)}{u}\;du
\eqno{(12)}$$

is defined for all $x \neq 0$.  (The Cauchy principal
value is taken in (12) if $x < 0$.)  However, the power series
$$
E_1(x) + \gamma + \ln x = \sum^\infty_{j=1}\;
   \frac{x^j(-1)^{j-1}}{j!j}
\eqno{(13)}$$

is unsatisfactory as a means of evaluating $E_1(x)$ for large positive $x$,
for the reasons discussed in Section~3 in connection with the power series
for $\exp(x)$.  For sufficiently large $x$ it is preferable
to use the asymptotic expansion [12]
$$ 
E_1(x) = \exp(-x)\; \sum^k_{j=1}\; \frac{(j-1)!(-1)^{j-1}}{x^j}
                           + R_k(x) \;,
\eqno{(14)}$$

where
$$ 
R_k(x) = k!(-1)^k\; \int^\infty_x \frac{\exp(-u)}{u^{k+1}}\;du \;.
\eqno{(15)}$$

For large positive $x$, the relative error attainable by using (14)
with $ k \simeq x$ is $O(x^{{1/2}} \exp(-x))$, because
$$
|R_k(k)|\;\leq\; k!k^{-(k+1)}\exp(-k)
            \; = \;O(k^{-{1/2}}\exp(-2k)) \;.
\eqno{(16)}$$

Thus, the asymptotic series may be used to evaluate $E_1(x)$ to precision
$n$ when\\
$x > n \ln 2 + O(\ln n)$.
(Similarly if $(-x) > n \ln 2 + O(\ln n)$,
although the estimation of $|R_k(-k)|$ is more
difficult than that of $|R_k(k)|$.)\\

There are many other examples where asymptotic expansions are useful, e.g.
for $\erf c(x)$ (mentioned in Section3), for Bessel functions [11, 35],
etc.  Asymptotic expansions aften arise when the convergence of series is
accelerated by the Euler-Maclaurin sum formula [1].  For example,
the Riemann zeta function $\zeta(s)$ is defined for $R(s) > 1$ by
$$ 
\zeta(s) = \sum^\infty_{j=1}\; j^{-s} \;,
\eqno{(17)}$$

and by analytic continuation for other $ s \neq 1$.  (Here we allow
complex $s$.)  $\zeta (s)$ may be evaluated to any desired precision if $m$
and $p$ are chosen large enough in the Euler-Maclaurin formula [8]
$$ 
\zeta(s) = \sum^{p-1}_{j=1}\; j^{-s} + {\ts \frac{1}{2}}p^{-s}
                + \frac{p^{1-s}}{s-1} +
                \sum^m_{k=1}\; T_{k,p}(s)+E_{m,p}(s) \;,
\eqno{(18)}$$
where
$$ T_{k,p}(s) = \frac{B_{2k}}{(2k)!}\; p^{1-s-2k} \prod^{2k-2}_{j=0}
                              \; (s+j) \;,
\eqno{(19)}$$

$$
|E_{m,p}(s)|\; <\;|T_{m+1,p}(s) \;
                 (s+2m+1)/(\sigma +2m+1)|\;,
\eqno{(20)}$$

$m \ge 0$, $p \ge 1$, $\sigma = R(s) > - (2m+1)$, and the $B_{2k}$ are
Bernoulli numbers.\\

In arbitrary-precision computations we must be able to compute as many terms 
of an asymptotic expansion as are required to give the desired accuracy.
It is easy to see that $m$ in (18) can not be bounded as the precision 
$n \rightarrow \infty$, else $p$ would have to increase as an
exponential function of $n$.  To evaluate $\zeta(s)$ from (18) to
precision $n$ in time polynomial in $n$, both $m$ and $p$ must
tend to infinity with $n$.  Thus, the Bernoulli numbers 
$B_2,\ldots,B_{2m}$ can not be stored in a table of fixed size,
but must be computed when needed (see Sections 7 and 9).  For this
reason we can not use asymptotic expansions when the general form
of the coefficients is unknown (such as Stirling's formula for
$\Gamma(x)$) in arbitrary-precision calculations.  Often there is a 
related expansion with known coefficients, e.g. the asymptotic
expansion for $\ln \Gamma(x)$ has coefficients related to the
Bernoulli numbers, like (19).
 
\section{Continued fractions}

Sometimes continued fractions are preferable to power series or
asymptotic expansions.  For example, Euler's continued fraction [34]
$$ 
\exp(x) \; E_1(x) = 1/x+1/1+1/x+2/1+2/x+3/1+\cdots
\eqno{(21)}$$

converges for all real $x>0$, and is better for computation of $E_1(x)$
than the power series (13) in the region where the power series
suffers from catastrophic cancellation but the asymptotic expansion (14) is
not sufficiently accurate.  Convergence of (21) is slow if $x$ is small,
so (21) is preferred for precision $n$ evaluation of $E_1(x)$ only
when $x \in (c_1 n,\; c_2 n)$, $c_1 \simeq 0.1$, $c_2 \simeq \ln 2$.\\

It is well known that continued fractions may be evaluated by either 
forward or backward recurrence relations.  Consider the finite continued
fraction
$$ y = a_1/b_1 + a_2/b_2 + \cdots + a_k/b_k \;.
\eqno{(22)}$$

The backward recurrence is $R_k = 1$, $R_{k-1} = b_k$,
$$ R_j = b_{j+1}\; R_{j+1} + a_{j+2}\; R_{j+2}\;
\hspace*{9mm}(j= k-2,\ldots,0) \;,
\eqno{(23)}$$

and $y = a_1 R_1/R_0$. The forward recurrence is
$P_0 = 0$, $P_1 = a_1$, $Q_0 = 1$, $Q_1 = b_1$,
$$
\left. \begin{array}{l}
       P_j = b_j\; P_{j-1} + a_j\; P_{j-2}\\[2pt]
       Q_j = b_j\; Q_{j-1} + a_j\; Q_{j-2} \\
       \end{array} \right\} \hspace*{9mm}(j=2, \ldots, k) \;,
\eqno{(24)}$$

and $y= P_k/Q_k$.\\

The advantage of evaluating an infinite continued fraction such as (21)
via the forward recurrence is that $k$ need not be chosen in advance;
we can stop when $|D_k|$ is sufficiently small, where
$$
D_k = \frac{P_k}{Q_k} - \frac{P_{k-1}}{Q_{k-1}} \;.
\eqno{(25)}$$

The disadvantage of the forward recurrence is that twice as many
arithmetic operations are required as for the backward recurrence
with the same value of $k$.  There is a simple solution to this
dilemma if we are working with variable-precision floating-point
arithmetic which is much more expensive than single-precision
floating-point.  We use the forward recurrence with single-precision
arithmetic (scaled to avoid overflow/underflow) to estimate $k$,
then use the backward recurrence with variable-precision arithmetic.  
One trick is needed:  to evaluate $D_k$ using scaled single-precision
we use the recurrence
$$
\left.\begin{array}{l}
      D_1 = a_1/b_1\;,\\[2pt]
      D_j = - a_j Q_{j-2} D_{j-1}/Q_j \hspace*{9mm}(j= 2,3,\ldots)\\
      \end{array} \right\}
\eqno{(26)}$$

which avoids the cancellation inherent in (25).\\

In recent versions of the MP package [7] we have used the continued
fraction (21) in the manner just described, and similar continued
fractions could well be used for the computation of other special
functions.  Since power series and asymptotic series are generally
easier to analyse and program than continued fractions, we have
avoided continued fractions except where they are clearly superior
to the other methods.

\section{Recurrence relations}

The evaluation of special functions by continued fractions is a special
case of their evaluation by recurrence relations.  For example, the
Bessel functions $J_\nu(x)$ satisfy the recurrence relation
$$
J_{\nu-1}(x) + J_{\nu+1}(x) = \frac{2\nu}{x} J_\nu(x)
\eqno{(27)}$$

which may be evaluated backwards (compare (23)), using a normalisation
condition such as
$$
J_0(x) + 2 \;\sum^\infty_{\nu=1}\; J_{2\nu}(x) = 1 \;.
\eqno{(28)}$$

This seems to be the most effective method in the region where Hankel's 
asymptotic expansion is insufficiently accurate but the power series
$$
J_\nu(x) = \left({\frac{x}{2}}\right)^\nu \;\sum^\infty_{j=0}\;
       \frac{(-x^2/4)^j}{j!\; \Gamma(\nu+j+1)}
\eqno{(29)}$$

suffers from catastrophic cancellation.  For details see [17].\\

In Section~5 the constants $C_k=B_{2k}/(2k)!$ were required, where
the $B_{2k}$ are Bernoulli numbers.  The $C_k$ are defined by the
generating function
$$
\sum^\infty_{k=0} \; C_k \; x^{2k} = \frac{x}{e^x -1}+ \frac{x}{2} \;.
\eqno{(30)}$$

Multiplying both sides by $e^x -1$ and equating coefficients gives 
the recurrence relation
$$
\frac{C_k}{1!}+\frac{C_{k-1}}{3!}+\cdots+\frac{C_1}{(2k-1)!}
        \;=\; \frac{k-\frac{1}{2}}{(2k+1)!} \;,
\eqno{(31)}$$

which has often been used to evaluate Bernoulli numbers [21].\\

Unfortunately, forward evaluation of the recurrence (31) is numerically
unstable:  using precision $n$ the relative error in the computed
$C_k$ is of order $4^k2^{-n}$.  We shall  not prove this, but
shall indicate why such behaviour is to be expected.  Consider the
``homogeneous'' recurrence
$$
\frac{\widehat{C}_k}{1!} + \frac{\widehat{C}_{k-1}}{3!} + \cdots
      + \frac{\widehat{C}_1}{(2k-1)!} \;= \;0 \hspace*{9mm} (k\ge 2)
\eqno{(32)}$$

with $\widehat{C}_1 = 1$, and let
$$
\widehat{G}(x) = \sum^\infty_{k=1}\; \widehat{C}_k\; x^{2k}
\eqno{(33)}$$

be the generating function for the $\widehat{C}_k$.  It is easy to
show that
$$
\widehat{G}(x) = \frac{x^3}{\sinh x} \;.
\eqno{(34)}$$

Thus $\widehat{G}(x)$ has poles at $\pm i\pi$ , and
$$
|\widehat{C}_k|\; \ge\; K\pi^{-2k}
\eqno{(35)}$$

for some $K>0$.  This suggests that  an error of order $2^{-n}$
in an early value of $C_j$ propagates to give an (absolute) error
of order $2^{-n}\pi^{-2k}$ in $C_k$ for large $k$.  Since
$|C_k|\sim (2\pi)^{-2k}$, this absolute error corresponds
to a relative error of order $2^{2k-n}=4^k2^{-n}$ in $C_k$.\\

Despite its numerical instability, use of (31) may give the $C_k$ 
to acceptable accuracy if they are only needed to generate coefficients 
in an Euler-Maclaurin expansion whose successive terms diminish by at
least a factor of 4.  If the $C_k$ or $B_{2k}$ are required to 
precision $n$, either (31) must be used with sufficient guard digits, 
or a more stable recurrence must be used.  If we multiply both sides
of (30) by $\sinh(x/2)/x$ and equate coefficients, we get the
recurrence
$$
C_k + \frac{C_{k-1}}{3!\;4} + \cdots + \frac{C_1}{(2k-1)!\;4^{k-1}}
           = \frac{2k}{(2k+1)!\;4^k}
\eqno{(36)}$$

If (36) is used to evaluate $C_k$, using precision $n$ arithmetic,
the error is only $O(k^22^{-n})$. Thus, this method is currently 
used in the MP package instead of a method based on (31).

\section{Newton's method}

Newton's method and related zero-finding methods may be used to
evaluate a function if we have an algorithm for evaluation of the
inverse function.  For example, applying Newton's method to
$f(x) = y-x^{-m}$ (where $y$ is regarded as constant) gives the
iteration
$$
x_{j+1} = x_j + x_j(1-x^m_j\; y)/m \;,
\eqno{(37)}$$

which converges (from a sufficiently good initial approximation)
to $y^{-1/m}$.  Note that (37) does not involve divisions except
by the small integer $m$.\\

Similarly, applying Newton's method to $f(x) = \exp(x) -y$
gives the iteration
$$
x_{j+1}= x_j + y \;\exp(-x_j) - 1 \;.
\eqno{(38)}$$

which converges to $\ln y$ if $x_0$ is a sufficiently good initial
approximation.\\

Newton's method generally has second order convergence, so we may 
start with low precision and approximately double it at each
iteration.  Thus, the work required is of the same order as the
work for the final iteration.  Applied to (37) with $m=1$ and $2$,
this argument shows that reciprocals and square roots can be found
to precision $n$ in $O(M(n))$ operations.  For futher details, and a
comparison of the efficiencies of various root-finding methods for
variable-precision computations, see [4, 5].

\section{Contour integration}

In this section we assume that facilities for variable-precision
complex arithmetic are available.  Let $f(z)$ be holomorphic in
the disc $|z|< R$, $R > 1$, and let the power series for
$f$ be
$$
f(z) = \sum^\infty_{j=0}\; a_j \;z^j	%
\eqno{(39)}$$

From Cauchy's theorem [18] we have
$$
a_j = \frac{1}{2\pi i}\; \int_C \;\frac{f(z)}{z^{j+1}}\;dz \;,
\eqno{(40)}$$

where $C$ is the unit circle.  The contour integral in (40) may be
approximated numerically by sums
$$
S_{j,k} = \frac{1}{k}\; \sum^{k-1}_{m=0}\; 
     f(e^{2\pi im/k})e^{-2\pi ijm/k} \;.
\eqno{(41)}$$

From Cauchy's theorem, %
provided $j < k$ and the contour $C$ is enlarged slightly to enclose the
$k$-th roots of unity, we have
$$
S_{j,k} - a_j = \frac{1}{2\pi i} \int_C 
    \frac{f(z)}{(z^k - 1)z^{j+1}}\; dz 
$$
\vspace*{-9mm}
\begin{flushright}
(42)
\end{flushright}
\vspace*{-4mm}
$$
\hspace*{10mm}= a_{j+k} + a_{j+2k} + \cdots \;,
$$

so $|S_{j,k} - a_j|= O((R - \delta)^{-(j+k)})$ as $ k \rightarrow
\infty$, for any $\delta > 0$.\\

For example, let
$$
f(z) = \frac{z}{e^z -1} + \frac{z}{2}
\eqno{(43)}$$

as in Section~7, so $a_{2j} = B_{2j}/(2_j)!$ and $R = 2\pi$.  Then
$$
S_{2j,k} - \frac{B_{2j}}{(2j)!} = \frac{B_{2j+k}}{(2j+k)!} +
           \frac{B_{2j+2k}}{(2j+2k)!} + \cdots \;,
\eqno{(44)}$$

so we can evaluate $B_{2j}$ with relative error $O((2\pi)^{-k})$
by evaluating $f(z)$ at $k$ points on the unit circle.  (By symmetry and
conjugacy only $k/4 + 1$ evaluations are required if $k$ is a multiple
of four.)  If $\exp(-2\pi ijm/k)$ is computed efficiently from 
$\exp(-2\pi i/k)$ in the obvious way, the time required to evaluate
$B_2,\ldots,B_{2j}$ to precision $n$ is $O(jnM(n))$, and the space
required is $O(n)$.  The recurrence relation method of Section~7 requires
time only $O(j^2n)$, 			%
but space $O(jn)$.  Thus, the method of contour
integration is recommended if space is more important than time.\\

For further discussion of the contour integration method, see [26].

\section{Special methods}

In this section we mention two of a large number of ``special'' methods
which are useful but less generally applicable than the methods of
Sections 3 to 9.  The first such method is the conversion of a power
series which suffers from catastrophic cancellation to one which is
better behaved numerically.  One example, (9), has already been given.
Another example occurs with
$$
E(x)  = \int^x_0 \frac{(1-e^{-u})}{u}\;du =
       \sum^\infty_{j=1}\;\frac{x^j(-1)^{j-1}}{j!\;j}
\eqno{(45)}$$

(a series encountered in Section~5).  Multiplying by $\exp(x)$ and
using some well known identities, we find
$$
\exp(x)E(x) = \sum^\infty_{j=1}\; H_j\; x^j/j! \;,
\eqno{(46)}$$
where
$$
H_j = \sum^j_{m=1}\; \frac{1}{m} \;.
\eqno{(47)}$$

If $x$ is large and positive, the series in (46) is much better behaved
numerically than the series in (45).  For an application where $E(x)$
was required to high precision with $x$ a positive integer, see [11].
At first sight it appears that, in this application, the summation to 
precision $n$ of $k$ terms in the series (45) requires $O(kn)$
operations, while (46) requires $\Omega(kM(n))$ operations.  However,
by a ``summation by parts'' trick described in [11], this can be
reduced to $O(kn)$ operations.\\

Our second ``special'' method is the evaluation of $\pi$ and elementary 
functions by the arithmetic-geometric mean (AGM) iteration.  It is
well known that the AGM can be used to compute elliptic integrals,
but perhaps less well known that it can also be used to compute 
$\pi$ and elementary functions, and gives the fastest known methods
when the precision $n$ is very large [4, 6].\\

The AGM of two positive numbers $a_0$ and $b_0$ is
$a = {\ds\lim_{j\rightarrow\infty}}\; a_j = {\ds\lim_{j\rightarrow\infty}}\; b_j$,
where
$$ 
a_{j+1} = \frac{a_j + b_j}{2}
\eqno{(48)}$$
and
$$
b_{j+1} = \sqrt{a_j b_j} \;.
\eqno{(49)}$$

There is no essential loss of generality in assuming that $a_0 = 1$ and
$b_0 = \cos\phi$.  Gauss [16] showed that $2a = \pi/K(\phi)$, where
$$
K(\phi) = \int^{\pi/2}_{0} \; (1-\sin^2\; \phi\; \sin^2 \theta)^{-{1/2}}
\; d\theta
\eqno{(50)}$$

is the complete elliptic integral of the first kind.  A simple proof is
given in [27].\\

The AGM iteration converges quadratically:  if $\varepsilon_j =
1-b_j/a_j$ then
$$
\varepsilon_{j+1} = 1-2(1-\varepsilon_j)^{{1/2}}/(2-\varepsilon_j)=
\varepsilon^2_j/8 + O(\varepsilon^3_j) \;.
\eqno{(51)}$$

Using the AGM and an identity of Lagrange, we get a family of
quadratically convergent algorithms for the computation of $\pi$.  The
simplest of these is:\\

\hspace*{18mm}\begin{tabular}{ll}
\multicolumn{2}{l}{$a\; := \; 1\;;$}\\
\multicolumn{2}{l}{$b\; := \; 1/\sqrt{2}\;;$}\\	%
\multicolumn{2}{l}{$t\; := \; 1/4\;;$}\\
\multicolumn{2}{l}{$j\; := \; 1\;;$}\\
\ul{repeat}&\\
           &$y\; := \; a\;;$\\
           &$a\; := \; (a+b)/2\;;$\\
           &$b\; := \; \sqrt{b\times y}\;;$\\	%
           &$t\; := \; t-j\times (a-y)^2\;;$\\	%
           &$j\; := \; 2 \times j$\\
\ul{until} &$(a-b)\; <\; {\rm{tolerance}}\;;$\\
\ul{return}&$a^2/t.$				%
\end{tabular}\\

After $k$ iterations the error $|a^2/t-\pi|$ is about
$8\pi\exp(-2^k\pi)$, e.g.~$k=5$ gives error less than $10^{-42}$.
For further details see [6, 31].\\

In [3, 4, 6] it is shown how the AGM may be used to compute the
elementary functions $\exp(x)$, $\ln(x)$, $\atan(x)$, $\sin(x)$
etc.~to precision $n$ in $O(M(n) \log n)$ operations. The
 factor ``$\log n$'' arises because $O(\log n)$ iterations of the AGM
are required.  It is important to note that the AGM iteration is not
self-correcting, so the trick of starting with low precision and
doubling it on each iteration (as used in Section~8) is  not
applicable.

\section{Summary}

Many ``classical'' methods may be adapted for use in variable-precision
computations; others are not readily adaptable.  Since the performance
criteria are different in variable-precision applications, the best
method may be one which is not well-suited to fixed-precision
computations.  For example, it might be numerically unstable, and
thus require the working precision to be increased.  The examples
given in Sections 3 to 10 above are intended to illustrate the main
ideas of variable-precision algorithms.

\section{Acknowledgement}

Christian Reinsch kindly suggested that the use of (36) would be faster 
than the method described in Section~6.11 of [7].

\addcontentsline{toc}{section}{References}

\section*{References}

\begin{tabular}{rp{408pt}}		%

[1]&
M.A.~Abramowitz and I.A.~Stegun (eds.),
{\em Handbook of Mathematical Functions with Formulas, Graphs,
and Mathematical Tables}, National Bureau of Standards, Washington,
D.C., 1964 (reprinted by Dover, New York, 1965).\\

[2]&
R.E.~Bank and D.J.~Rose,
``Extrapolated fast direct algorithms for elliptic boundary value
problems'',
in {\em Algorithms and Complexity\/} (ed.~by J.F.~Traub), Academic Press,
New York, 1976, 201--249.\\

[3]&
M.~Beeler, R.W.~Gosper and R.~Schroeppel,
{\em Hakmem}, M.I.T. Artificial Intelligence Lab. Memo No.~239,
Feb.~1972.
[Available from\\
&{\tt http://{\lbrk}www.inwap.com/{\lbrk}pdp10/%
{\lbrk}hbaker/{\lbrk}hakmem/{\lbrk}hakmem.html}~]\\

[4]&
R.P.~Brent,							%
``Multiple-precision zero-finding methods and the complexity of
elementary function evaluation'',
in {\em Analytic Computational Complexity\/} (ed.~by J.F.~Traub),
Academic Press, New York, 1975, 151--176.
[Available from\\
&{\tt http://wwwmaths.anu.edu.au/{\realtilde}brent/pub/pub028.html}~]\\

[5]&
R.P.~Brent,							%
``The complexity of multiple-precision arithmetic'',
in {\em The Complexity of Computational Problem Solving\/}
(eds.~R.S.~Anderssen and R.P.~Brent),
Queensland Univ.~Press, Brisbane, 1976, 126--165.
[Available from\\
&{\tt http://wwwmaths.anu.edu.au/{\realtilde}brent/pub/pub032.html}~]\\

[6]&								%
R.P.~Brent,
``Fast multiple-precision evaluation of elementary functions'',
{\em  J.ACM\/} 23 (1976), 242--251.\\

[7]&								%
{R.P.~Brent, ``A Fortran multiple-precision arithmetic package'',}
{\em ACM Trans.\ Math.\ Software\/} 4 (1978), 57--70.\\

[8]&
R.P.~Brent,							%
``On the zeros of the Riemann zeta function in the critical strip'',
{\em Math.~Comp.} 33 (1979), 1361--1372.\\

[9]&
R.P.~Brent,
{\em Numerical investigation of the Riemann-Siegel approximation},
unpublished notes, 1979.\\			%

[10]&
R.P.~Brent, J.A.~Hooper and J.M.~Yohe,				%
``An Augment interface for Brent's multiple-precision arithmetic
package'',
{\em ACM Trans.\ Math.\ Software}, to appear.
[Appeared in vol.\ 6 (1980), 146--149.]\\

[11]&
R.P.~Brent and E.M.~McMillan,					%
``Some new algorithms for high-precision computation of Euler's
constant'',
{\em Math.~Comp.,} to appear.
[Appeared in vol.\ 34 (1980), 305--312.]\\

[12]&
N.G.~de Bruijn,
{\em Asymptotic Methods in Analysis},
3rd edition, North-Holland, 1970.\\

[13]&
C.W.~Clenshaw and F.W.J.~Olver,
``An unrestricted algorithm for the exponential function'',
{\em SIAM J.~Numer.~Anal.,} to appear.
[Appeared in vol.\ 17, 1980, 310--331.]\\

[14]&
J.T.~Coonen, W.~Kahan, J.~Palmer, T.~Pittman and D.~Stevenson,
``A proposed standard for binary floating point arithmetic,
draft 5.11'',
{\em ACM SIGNUM Newsletter}, October 1979, 4--12.\\

[15]&
G.E.~Forsythe,
``Pitfalls in computation, or why a math book isn't enough'',
{\em Amer. Math.~Monthly\/} 77 (1970), 931--956.\\

[16]&
C.F.~Gauss,
{\em Carl Friedrich Gauss Werke,}
Bd.~3, G\"{o}ttingen, 1876, 362--403.\\

[17]&
W.~Gautschi,
``Algorithm 236:  Bessel functions of the first kind'',
{\em Comm.~ACM\/} 7 (1964), 479--480.\\

[18]&
E.~Hille,
{\em Analytic Function Theory,} Vol.~1, Blaisdell, New York,
1959, Ch.~7.\\

[19]&
T.E.~Hull and J.J.~Hofbauer,
{\em Language facilities for multiple-precision floating-point
computation},
Dept.~of Computer Science, Univ.~of Toronto, 1974.\\

[20]&
J.P.~Jeter,
{\em A Variable-Precision Interval Data Type Extension to Fortran},
M.Sc.~thesis, Dept.~of Computer Science, Univ.~of S.W.~Louisiana,
Lafayette, Louisiana, July 1979.\\

[21]&
D.E.~Knuth,
``Euler's constant to 1271 places'',
{\em Math.~Comp.\/} 16 (1962), 275--281.\\

\end{tabular}

\begin{tabular}{rp{408pt}}

[22]&
D.E.~Knuth,
{\em The Art of Computer Programming,} Vol.~2,
Addison Wesley, Reading, Mass., 1969.\\

[23]&
D.E.~Knuth,
``Big Omicron and Big Omega and Big Theta'',
{\em SIGACT News\/} 8, 2 (1976), 18--24.\\

[24]&
D.H.~Lehmer,
``Tables to many places of decimals'',
{\em Math.~Tables Aids Comput.\/} 1 (1943), 30--31.
[The journal is now called {\em Mathematics of Computation}.]\\

[25]&
Y.L.~Luke,
{\em Algorithms for the Computation of Mathematical Functions},
Academic Press, New York, 1977.\\

[26]&
J.N.~Lyness and C.B~Moler,
``Numerical differentiation of analytic functions'',
{\em SIAM J.~Numer.~Anal.} 4 (1967), 202--210.\\

[27]&
Z.A.~Melzak,
{\em Companion to Concrete Mathematics},
Wiley, New York, 1973.\\

[28]&
R.E.~Moore,
{\em Interval Analysis},
Prentice-Hall, New Jersey, 1966.\\

[29]&
M.~Payne and W.~Strecker,
``Draft proposal for a binary normalized floating point standard'',
{\em ACM SIGNUM Newsletter}, October 1979, 24--30.\\

[30]&
H.J.J.~te Riele,
``Computations concerning the conjecture of Mertens'',
{\em J.~reine angew.~Math.} 311/312 (1979), 356--360.\\

[31]&
E.~Salamin, 
``Computation of $\pi$ using arithmetic-geometric mean'',
{\em  Math.~Comp.} 30 (1976), 565--570.\\

[32]&
J.L.~Schonfelder and J.T.~Thomason,
{\em Applications support by direct language extension -- an arbitrary
precision arithmetic facility in Algol 68},
Comp.~Centre, Univ.~of Birmingham, Birmingham, U.K., 1975.\\

[33]&
A.~Sch\"{o}nhage and V.~Strassen, 
``Schnelle Multiplikation grosser Zahlen'',
{\em Computing\/} 7 (1971), 281--292.\\

[34]&
H.~Wall,
{\em Analytic Theory of Continued Fractions\/},
Van Nostrand, New York, 1948.\\

[35]&
E.T.~Whittaker and G.N.~Watson,
{\em A Course of Modern Analysis\/},
Cambridge Univ.~Press, 1902.\\

[36]&
J.H.~Wilkinson,
{\em Rounding Errors in Algebraic Processes\/},
HMSO, London, 1963.\\

[37]&
J.M.~Yohe,
{\em The interval arithmetic package -- multiple precision version\/},
MRC Tech. Summary Report No.~1908, Math.~Res.~Center, Madison, Jan.~1979.
\end{tabular}

\end{document}